\documentclass[conference]{IEEEtran}

\IEEEoverridecommandlockouts                              

\usepackage{authblk}

\usepackage{cite}
\usepackage{comment}
\usepackage{subcaption}
\usepackage{epsfig} 
\usepackage{amsmath} 
\usepackage{amssymb}  
\usepackage{tikz}
\usetikzlibrary{shapes.geometric, arrows, backgrounds, scopes, positioning}
\usepackage{fourier}
\usepackage{mdframed}
\usepackage{algorithm}
\usepackage{algpseudocode}
\usetikzlibrary{fit}
\usepackage{hyperref}

\providecommand{\norm}[1]{\ensuremath{\left\lVert#1\right\rVert }}
\newtheorem{theorem}{\textbf{Theorem}}
\newtheorem{lemma}{\textbf{Lemma}}

\newtheorem{proposition}{\textbf{Proposition}}

\def\R{\mathbb{R}}

\def\A{A^TA}

\def\ro{\rho^*_{\beta}}
\DeclareRobustCommand{\bigO}{%
  \text{\usefont{OMS}{cmsy}{m}{n}O}%
}

\title{
Iterative Pre-Conditioning to Expedite the Gradient-Descent Method
}

\author{Kushal Chakrabarti$^\star$, Nirupam Gupta$^\dagger$, and Nikhil Chopra$^\star$
\thanks{$^\star$ University of Maryland, College Park, Maryland 20742, U.S.A. \\
$^\dagger$ Georgetown University, Washington, DC 20057, U.S.A. \\
Emails: {\em kchakrabarti0@gmail.com}, {\em nirupam.gupta@georgetown.edu} and {\em nchopra@umd.edu}}%
}

\begin{document}

\maketitle
\thispagestyle{empty}
\pagestyle{empty}

\begin{abstract}
This paper considers the problem of multi-agent distributed optimization. In this problem, there are multiple agents in the system, and each agent only knows its local cost function. The objective for the agents is to collectively compute a common minimum of the aggregate of all their local cost functions. In principle, this problem is solvable using a distributed variant of the traditional gradient-descent method, which is an iterative method. However, the speed of convergence of the traditional gradient-descent method is highly influenced by the {\em conditioning} of the optimization problem being solved. Specifically, the method requires a large number of iterations to converge to a solution if the optimization problem is {\em ill-conditioned}. 

In this paper, we propose an {\em iterative pre-conditioning} approach that can significantly attenuate the influence of the problem's conditioning on the convergence-speed of the gradient-descent method. The proposed pre-conditioning approach can be easily implemented in distributed systems and has minimal computation and communication overhead. For now, we only consider a specific distributed optimization problem wherein the individual local cost functions of the agents are {\em quadratic}. Besides the theoretical guarantees, the improved convergence speed of our approach is demonstrated through experiments on a real data-set.


\end{abstract}

\section{Introduction}


We consider a synchronous distributed system that comprises of a server and $m$ agents in a server-based architecture, as shown in Fig.~\ref{fig:sys}. The server-based architecture can be emulated easily on a rooted peer-to-peer network using the well-known routing and message authentication primitives~\cite{lynch1996distributed, tanenbaum1981network}. Each agent $i \in \{1, \ldots, \, m\}$ holds a pair $(A^i, \, b^i)$, where $A^i$ is a matrix of size $n_i \times n$ and $b^i$ is a column-vector of size $n$. Let $\R^n$ denote the set of real-valued vectors of size $n$. For a vector $v \in \R^n$, let $\norm{v}$ denote its 2-norm. If $(\cdot)^T$ denotes the transpose then $\norm{v}^2 = v^T v$. The objective of the agents is to solve for the least-squares problem:
\begin{align}
\operatorname*{\textit{minimize}}_{x\in \R^n} ~ & \sum_{i=1}^m \dfrac{1}{2}\norm{A^i \, x - b^i}^2. \label{eqn:opt_1}
\end{align}
Applications of such problems include linear regression, state estimation, hypothesis testing.

\tikzstyle{server} = [rectangle, rounded corners, minimum width=1.7cm, minimum height=1cm,text centered, text width=1cm, draw=black, fill=blue!30]
\tikzstyle{agent} = [rectangle, minimum width=1cm, minimum height=1cm,text centered, text width=1.8cm, draw=black, fill=blue!10]
\tikzstyle{dots} = [circle, inner sep=0pt,minimum size=2pt, draw=black, fill=blue!50!cyan]
\tikzstyle{arrow} = [thick,<->,>=stealth]

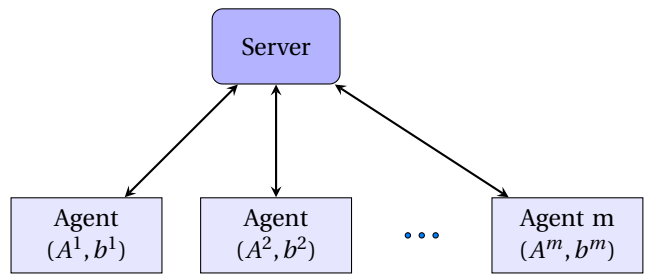
\begin{figure}[thpb]  
\centering
\begin{tikzpicture}[node distance = 1.5cm, auto]
    \node (server) [server] {Server};
    \node (m/c2) [agent, below  = of server] {Agent\\ $(A^2,b^2)$};
    \node (m/c1) [agent, left = of m/c2, xshift=1cm] {Agent\\ $(A^1,b^1)$};
    \node (d1) [dots, right = of m/c2, xshift=-0.8cm] {};
    \node (d2) [dots, right = of d1, xshift=-1.4cm] {};
    \node (d3) [dots, right = of d2, xshift=-1.4cm] {};
    \node (m/c3) [agent, right = of d3, xshift=-0.8cm] {Agent m\\ $(A^m,b^m)$};
    
    \draw[arrow] (server) -- (m/c1);
    \draw[arrow] (server) -- (m/c2);
    \draw[arrow] (server) -- (m/c3);
\end{tikzpicture}
\caption{The system architecture.}
\label{fig:sys}
\end{figure}

The server collaborates with the agents to solve the problem. The server-based architecture is also commonly known as the {\em federated model}~\cite{yang2019federated}, where the agents collaboratively solve an optimization problem, such as~\eqref{eqn:opt_1}, using a central server without ever sharing their local data $(A^i, \, b^i)$ with anyone in the system. In this architecture, the server maintains an estimate of the global solution of~\eqref{eqn:opt_1}, and each agent downloads the estimated global solution from the server. Each agent updates the estimated global solution locally using the data it possesses and uploads its local updated estimate to the server. The server accumulates the locally updated solutions from all the agents to improve the global solution. \\


The most common and most straightforward algorithm for the above distributed minimization problem is the distributed variant of the vanilla gradient-descent method, called the distributed gradient-descent (DGD) method~\cite{bertsekas1989parallel}. Recently, Azizan-Ruhi et al., 2019~\cite{azizan2019distributed} have proposed an {\em accelerated} projection method, which is built upon the seminal work on accelerated methods by Nesterov~\cite{nesterov27method}. However, Azizan-Ruhi et al.~do not provide any theoretical guarantees on the improved convergence of their method over the traditional DGD method~\cite{azizan2019distributed}. \\

In this paper, we propose a different approach than the accelerated methods to improve the convergence speed of the DGD method. Instead of using {\em momentum}, we use a {\em pre-conditioning} scheme wherein each iteration the update vector, used by the server to update its estimate of the global solution of~\eqref{eqn:opt_1} is multiplied by a {\em pre-conditioner matrix} that we provide. The pre-conditioner matrix in itself is updated in each iteration; hence we call our approach the {\em iterative pre-conditioning} method.\\

The DGD method is an iterative algorithm, in which the server maintains an estimate for a point of minimum of~\eqref{eqn:opt_1}, denoted by $x^{*}$, and updates it iteratively by collaborating with the agents as follows. For each iteration $t = 0, 1,\, \dots$, let $x(t)$ denote the estimate of $x^{*}$ at the beginning of iteration $t$. The initial value $x(0)$ is chosen arbitrarily. In each iteration $t$, the server broadcasts $x(t)$ to all the agents. Each agent $i$ computes the gradient of the function $(1/2) \norm{A^i x - b^i}^2$ at $x = x(t)$, denoted by $g^i(t)$, and sends it to the server. Note that
\begin{align}
    g^i(t) = (A^i)^T \, \left(A^i \, x(t) - b^i\right), \quad \forall i, \, t. \label{eqn:g_i}
\end{align}
Upon receiving the gradients $\{g^i(t) \, | \, i = 1, \ldots, \, m\}$ from all the agents, the server updates $x(t)$ to $x(t+1)$ using step-size of constant value $\delta$ as follows:
\begin{align}
    x(t+1) = x(t) - \delta \, \sum_{i = 1}^m g^i(t), \quad \forall t. \label{eqn:dgd}
\end{align}
To be able to present the contribution of this paper, we first briefly review the convergence of the DGD method described above.\\

We define the following notations. Let $A = [(A^1)^T, \ldots, \, (A^m)^T]^T $ denote the matrix obtained by stacking matrices $A^i$'s vertically. So, matrix $A$ is of size $N \times n$ where $N = (\sum_{i = 1}^m n_i)$. We assume that $A$ is a tall matrix, i.e., $N \geq n$. Similarly, we concatenates the $b^i$'s to get $b = \begin{bmatrix} (b^1)^T, \ldots, \, (b^m)^T \end{bmatrix}^T \in \R^N$.

\subsection{Convergence of DGD}
\label{sub:conv_dgd}

If matrix $A$ is full-column rank then we know that there exists step-size $\delta$ for which there is a positive value $\rho < 1$ such that~\cite{fessler2008image},
\begin{align*}
    \norm{x(t) - x^*} \leq \rho^{t} \norm{x(0) - x^*}, \quad t = 0, \, 1, \ldots
\end{align*}
The value $\rho$ is commonly referred as the {\em convergence rate}. Smaller $\rho$ implies higher {\em convergence speed}, and vice-versa. If we let $\lambda$ and $\gamma$ denote the largest and smallest eigenvalues of $A^T A$, then it is known that
\begin{align}
    \rho \geq \rho_{GD} = \frac{(\lambda/\gamma) - 1}{(\lambda/\gamma) + 1}. \label{eqn:conv_dgd}
\end{align}
The ratio $\lambda/\gamma$ is also commonly referred as the \emph{condition number} of matrix $A^T A$ which we denote by $\kappa(A^T A)$. 

\subsection{Pre-Conditioning of DGD}
\label{sub:pre_dgd}
Our objective is to improve the convergence rate $\rho$ of DGD beyond $\rho_{GD}$ by a suitable {\em pre-conditioning} proposed in this paper. Let $K$, referred as the {\em pre-conditioner}, be a square matrix of size $n \times n$. The server now updates its estimate as follows:
\begin{align}
    x(t+1) = x(t) - \delta \, K \, \sum_{i = 1}^m g^i(t), \quad \forall t. \label{eqn:pre-con_dgd}
\end{align}
If the matrix product $K \A$ is positive definite then convergence of~\eqref{eqn:pre-con_dgd} can be made linear by choosing $\delta$ appropriately, and the smallest possible convergence rate for~\eqref{eqn:pre-con_dgd} is given by (ref. Chapter 11.3.3 of~\cite{fessler2008image})
\begin{align} \label{eqn:rho_K}
    \rho^*_{K} = \frac{\kappa(K \, A^T A) - 1}{\kappa(K \, A^T A) + 1}.    
\end{align}
However, most of the existing pre-conditioning techniques are not applicable for the distributed framework considered in this paper. The incomplete LU factorization algorithms~\cite{meijerink1977iterative}, accelerated iterative methods~\cite{axelsson1985survey}, symmetric successive over-relaxation method~\cite{axelsson1985survey},
for computing such a matrix $K$ such that $\kappa(K \, A^T A)$ is provably smaller than $\kappa(A^T A)$ require the server to have direct access to the matrices $A^i$'s. 
Some other pre-conditioning methods~\cite{benzi2002preconditioning} require $A$ to be a symmetric positive definite matrix. The recently proposed distributed pre-conditioning scheme for the D-Heavy Ball method in~\cite{azizan2019distributed} has the same convergence rate as APC. The experimental results in Section~\ref{sec:exp} suggest that our proposed scheme converges faster than APC, and hence faster than the said pre-conditioning scheme.

\subsection{Summary of Our Contributions}
\label{sub:contri}
We propose an \emph{iterative pre-conditioner} matrix $K(t)$, instead of a constant pre-conditioner matrix $K$. That is, the server updates its estimate as follows:
\begin{align}
    x(t+1) = x(t) - \delta \, K(t) \, \sum_{i = 1}^m g^i(t), \quad \forall t. \label{eqn:prop_dgd}
\end{align}
The pre-conditioner $K(t)$ can be computed in a distributed manner in the federated architecture for each iteration $t$, as is presented in Section~\ref{sec:algo}. We show that the iterative process~\eqref{eqn:prop_dgd} converges provably faster to the optimum point $x^*$ than the original DGD algorithm~\eqref{eqn:dgd}. In the experiments, we have also observed that the convergence speed of our proposed algorithm is faster than the accelerated projection-based consensus (APC) method proposed in~\cite{azizan2019distributed}. \\

Moreover, the computational complexity of the proposed method is the same as the DGD algorithm, which is the smallest amongst the existing distributed algorithms for solving~\eqref{eqn:opt_1}. We have formally shown the proposed algorithm to converge faster than the distributed gradient method. In contrast, the conditions where APC is guaranteed to converge faster than DGD have not been provided in~\cite{azizan2019distributed}. APC, as well as the pre-conditioning scheme for D-HBM in~\cite{azizan2019distributed}, have additional computational overhead before the iterations begin.\\

\section{PROPOSED ALGORITHM}
\label{sec:algo}

In this section, we present our proposed algorithm and theoretically guarantee that the proposed method converges faster than the distributed gradient descent method~\eqref{eqn:dgd} for solving the distributed least-squares problem~\eqref{eqn:opt_1}. \\

To be able to present our algorithm, we introduce the following notation.
\begin{itemize}
\setlength\itemsep{0.5em}
    \item For a positive integer $m$, let 
    \[[m]:= \{1,\ldots,m\}.\]
    \item Let $e_j$ denote the $j$-th column of the $(n \times n)$-dimensional identity matrix $I$.
    \item Let $K(t) \in \R^{n \times n}$ be the pre-conditioner matrix for iteration $t$, and let $k_j(t) \in \R^n$ denote the $j$-th column of $K(t)$.
    \item For $\beta > 0$, define
    \begin{align}\label{eqn:Rij}
        R^i_j(t) \hspace{-0.2em} & := \hspace{-0.2em}\left((A^i)^T \hspace{-0.2em}A^i + \dfrac{\beta}{m}I\right) k_j(t) \hspace{-0.1em} -  \dfrac{1}{m} e_j, \, j=1,...,n, \, i \in [m].
    \end{align}
\end{itemize}

The proposed method is summarized in Algorithm~\ref{algo_1}. Note that $\alpha, \,\delta, \, \beta$ are positive valued parameters in the algorithm.\\
\begin{algorithm}
  \caption{}\label{algo_1}
  \begin{algorithmic}[1]
    \State Initialize $x(0) \in \R^n$, $K(-1) \in \R^{n \times n}$
    \For{\texttt{$t=0,1,...$}}
      \State The server transmits $x(t)$ and $K(t-1)$ to all the agents $i \in [m]$
      \For{each agent $i \in [m]$}
         \State Compute $g^i(t)$ using~\eqref{eqn:g_i}
         \For{each $j \in {1,...,n}$}
           \State Compute $R^i_j(t-1)$ given by~\eqref{eqn:Rij}
         \EndFor
         \State Transmit $g^i(t)$ and $\{R^i_j(t-1)\}_{j=1}^n$ to the server
      \EndFor
      \For{each $j \in {1,...,n}$}
         \State The server updates each column of $K(t)$
         \begin{align}
             k_j(t) = k_j(t-1) - \alpha \sum_{i=1}^m R^i_j(t-1) \label{eqn:kcol_update}
         \end{align}
      \EndFor
      \State The server updates the estimate 
      \begin{align}
             x(t+1) = x(t) - \delta K(t) \sum_{i=1}^m g^i(t) \label{eqn:x_update}
         \end{align}
    \EndFor
  \end{algorithmic}
\end{algorithm}

For each iteration $t = 0, \, 1, \ldots$, the server maintains a pre-conditioner matrix $K(t-1)$ and an estimate $x(t)$ of the point of optimum $x^*$. The initial pre-conditioner matrix $K(-1)$ and estimate $x(0)$ are chosen arbitrarily. The server sends $x(t)$ and $K(t-1)$ to the agents. Each agent $i \in [m]$ computes $R^i_j(t-1)$ for $j = 1, \ldots, \, n$ as given in~\eqref{eqn:Rij}, the local gradient $g^i(t)$ as given in~\eqref{eqn:g_i}, and sends these to the server. Then, the server computes the updated pre-conditioner $K(t)$ as given in~\eqref{eqn:kcol_update}, and uses this updated pre-conditioner $K(t)$ to compute the updated estimate $x(t+1)$ as given by~\eqref{eqn:x_update}. \\

\subsection{Computational Complexity}

We count the number of floating-point multiplications needed per-iteration of Algorithm~\ref{algo_1}. 
For each agent $i$, the computation of $\{\R^i_j(t) | j=1,...,n\}$ in~\eqref{eqn:Rij} are independent of each other. So, each agent $i$ computes them in parallel. Now,~\eqref{eqn:Rij} can be rewritten as
\begin{align*}
        R^i_j(t) \hspace{-0.2em} & := \hspace{-0.2em}(A^i)^T \hspace{-0.2em}A^i k_j(t) + \dfrac{\beta}{m} k_j(t) \hspace{-0.1em} -  \dfrac{1}{m} e_j.
    \end{align*}
Two matrix-vector multiplications are required for computing $R^i_j(t)$: $A^i \, k_j(t)$ and $(A^i)^T(A^i k_j(t))$, in that order. Here $A$ is an $(n_i \times n)$-dimensional matrix and $k_j(t)$, $A^i \, k_j(t)$ respectively are column vectors of dimension $n$, $N$. So, a total of $(2n_i n+n)$ flops are required for computing $R_j(t)$. 

Two matrix-vector multiplications are required for computing $g^i(t)$: $A^i \, x(t)$, $(A^i)^T \, \left(A^i \, x(t) - b^i\right)$, in that order. By a similar argument as for $R^i_j(t)$, computing $g^i(t)$ needs $\bigO(n_in)$ flops.
One matrix-vector multiplication is required for computing $x(t+1)$: $K(t) \sum_{i=1}^m g^i(t)$, which needs $\bigO(n_in)$ flops. So, per-iteration computational cost of Algorithm~\ref{algo_1} is $\bigO(n_in)$.


\subsection{Convergence Analysis}
To be able to present the convergence of Algorithm~\ref{algo_1}, we introduce the following notations.
\begin{itemize}
\setlength\itemsep{0.5em}
    \item Let, \[K^* = \left(\A+\beta I\right)^{-1}.\] 
    The matrix $A^T A$ is positive semi-definite. Thus, $\left(\A+\beta I\right)$ is positive definite for $\beta > 0$. So, $K^*$ is well-defined.
    \item Let $k_j^*$ be the $j$-th column of $K^*$, $j=1,\ldots,n$.
    \item Let, $\lambda$ and $\gamma$ denote the smallest and largest eigenvalue of $\A$, respectively.
    \item Let,
    \begin{align*}
        \rho^*_K & := \frac{\lambda - \gamma }{\lambda + \gamma +2\beta}, \, \ro := \dfrac{(\lambda - \gamma)\beta }{(\lambda + \gamma)\beta + 2\lambda \gamma}, \\ 
        \sigma_0 & := \delta \lambda \norm{K(-1)-K^*}_F.
    \end{align*}
\end{itemize}

We make the following assumption.\\
\noindent \textbf{Assumption 1}: Assume that the matrix $\A$ is full rank.

The above lemma shows that each column of the pre-conditioner matrix $K(t)$ asymptotically converges to the corresponding column of $K^*$. In other words, the matrix $K(t)$ asymptotically converges to $K^*$.\\

\begin{lemma} \label{thm:lem1}
Consider the iterative process~\eqref{eqn:kcol_update}, and let $\beta > 0$. Under Assumption 1, there exists $\alpha > 0$ for which there is a positive value $\rho_K < 1$ such that for each $j=1,\ldots,n$,
\begin{align*}
    \norm{k_j(t)-k_j^*} \leq \rho_K \norm{k_j(t-1)-k_j^*}, \, t=0,1,2,\ldots,
\end{align*}
where $\rho_K \geq \rho^*_K$.
\end{lemma}


The smallest and the largest Eigenvalues of $K^*\A$ have a direct influence on the convergence speed of Algorithm~\ref{algo_1}. The following lemma finds these Eigenvalues.\\

\begin{lemma} \label{thm:lem2}
If Assumption 1 holds and $\beta > 0$, then $K^*\A$ is positive definite, and the largest and the smallest eigenvalues of $K^*\A$ are $\dfrac{\lambda}{\lambda+\beta}$ and $\dfrac{\gamma}{\gamma+\beta}$, respectively.
\end{lemma}


The following proposition tells us that Algorithm~\ref{algo_1} converges to the point of minima at a linear rate. \\

\begin{proposition} \label{thm:pro1}
Consider Algorithm 1 with $\beta > 0$. Under Assumption 1, there exists $\delta > 0$ for which there is a positive value $\rho_{\beta} < 1$ such that
\begin{align*}
    \norm{x(t+1)-x^*} \leq \left(\rho_{\beta}+\sigma_0 (\rho^*_K)^{t+1}\right) \norm{x(t)-x^*}, \, \forall t,
\end{align*}
where $\rho_{\beta} \geq \ro$.
\end{proposition}


From Proposition \ref{thm:pro1}, the rate of convergence of Algorithm~\ref{algo_1} depends on the choice of initial value $K(-1)$ through $\sigma_0$. The closer $K(-1)$ is to $K^*$, the faster is the rate. \\

\begin{figure*}[htb!]
\begin{subfigure}{.45\textwidth}
  \centering
  \includegraphics[width = \textwidth]{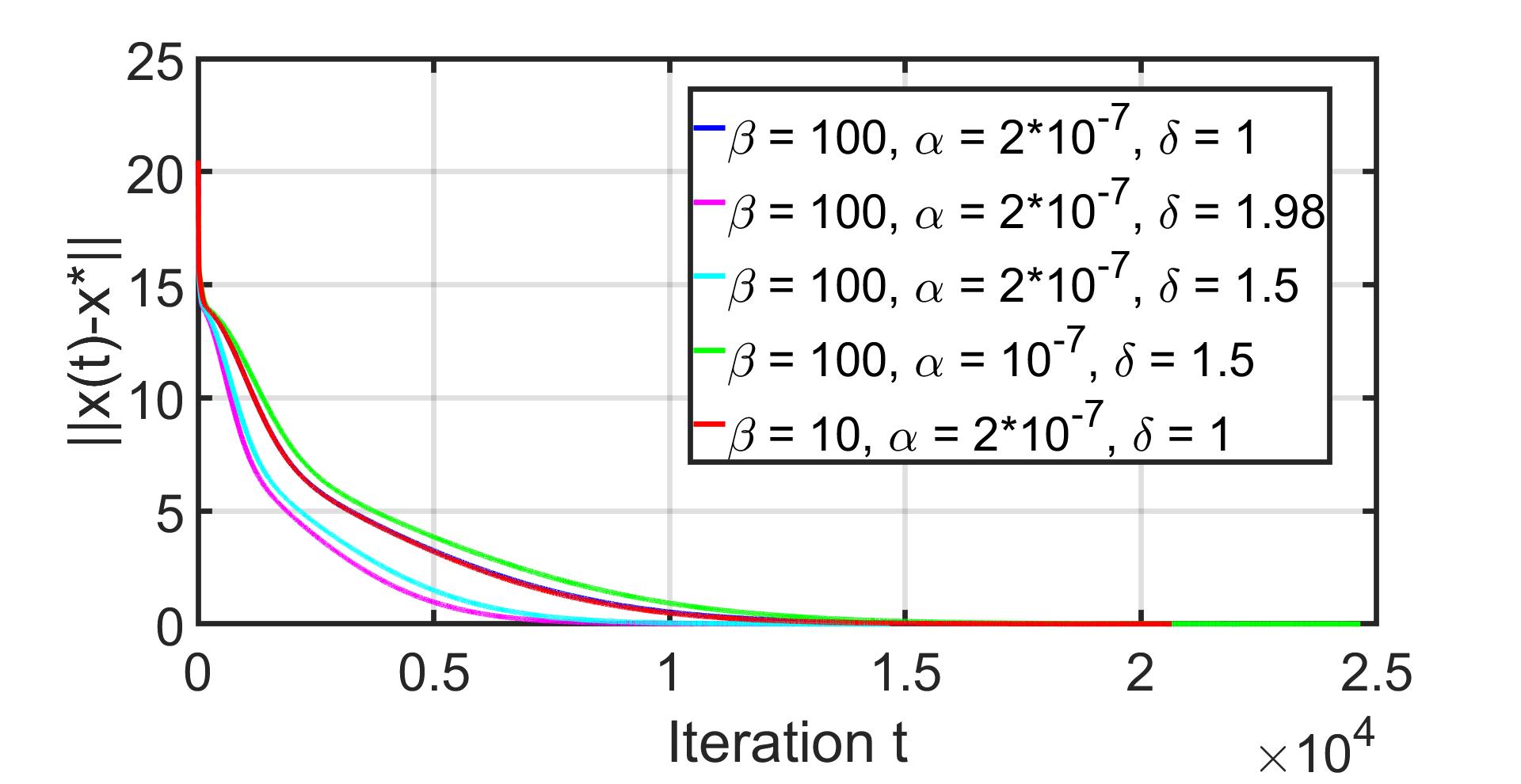}
  \caption{}
  \label{fig:bcsstm_1a}
\end{subfigure}%
\begin{subfigure}{.45\textwidth}
  \centering
  \includegraphics[width = \textwidth]{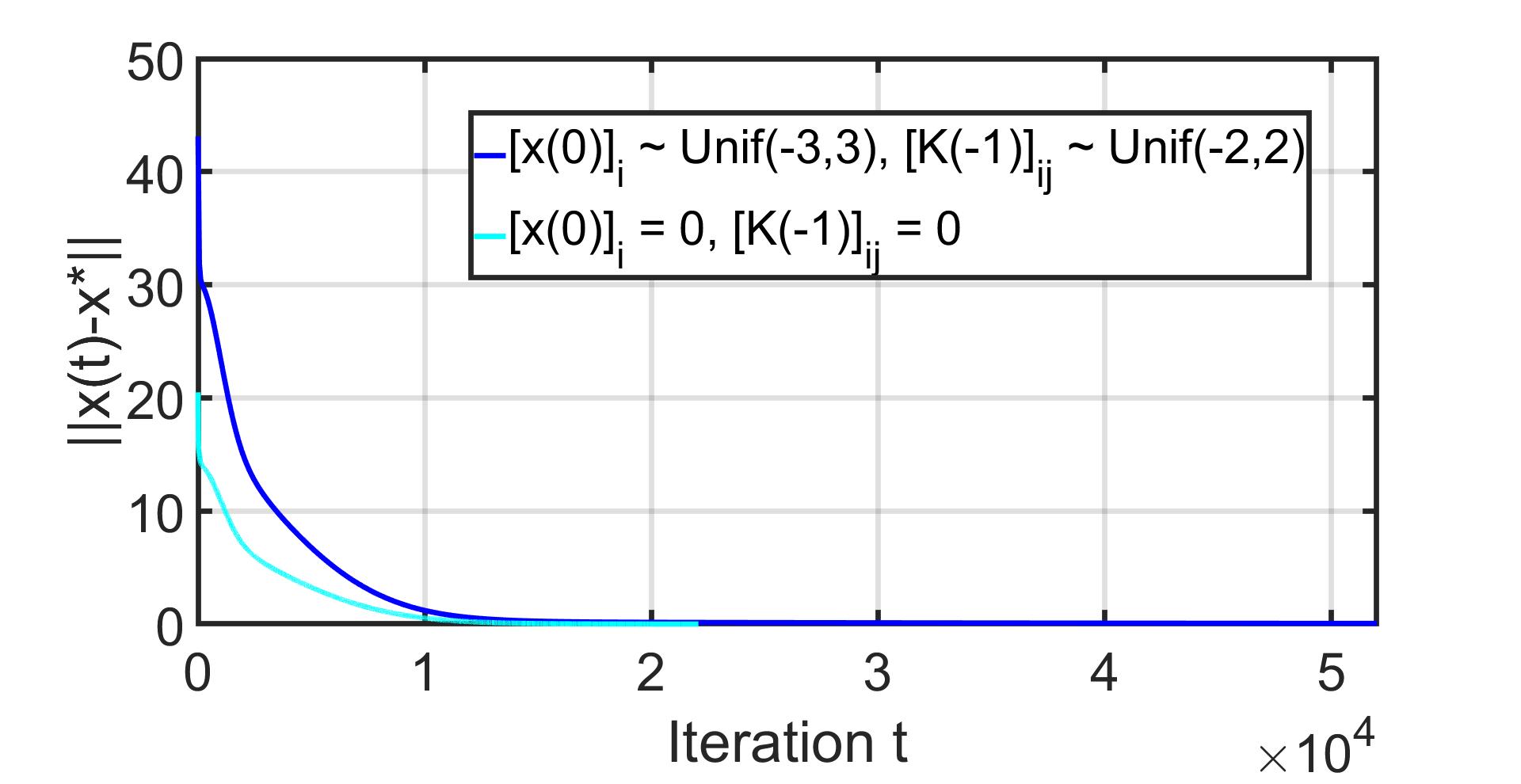}
  \caption{}
  \label{fig:bcsstm_1b}
\end{subfigure}
\caption{\footnotesize{Temporal evolution of error norm for estimate $\norm{x(t)-x^*}$ for Example 1, under Algorithm 1 with different parameter choices and initialization. (a) $x(0) = [0,\ldots,0]^T$; (b) $\beta = 100, \alpha = 2*10^{-7}, \delta = 1$.}}
\label{fig:bcsstm07}
\end{figure*}

\begin{figure*}[htb!]
\begin{subfigure}{.45\textwidth}
  \centering
  \includegraphics[width = \textwidth]{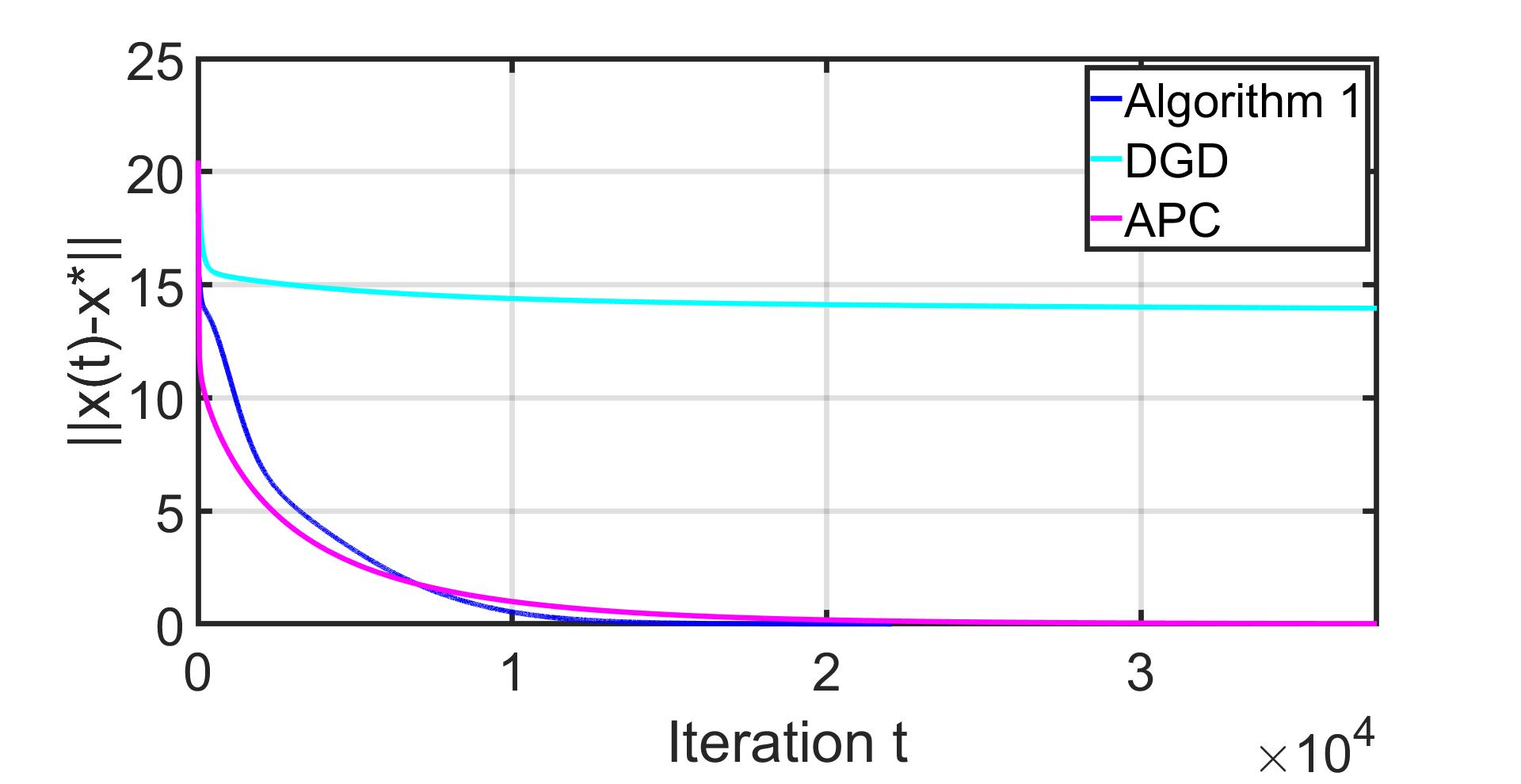}
  \caption{}
  \label{fig:comp1}
\end{subfigure}%
\begin{subfigure}{.45\textwidth}
  \centering
  \includegraphics[width = \textwidth]{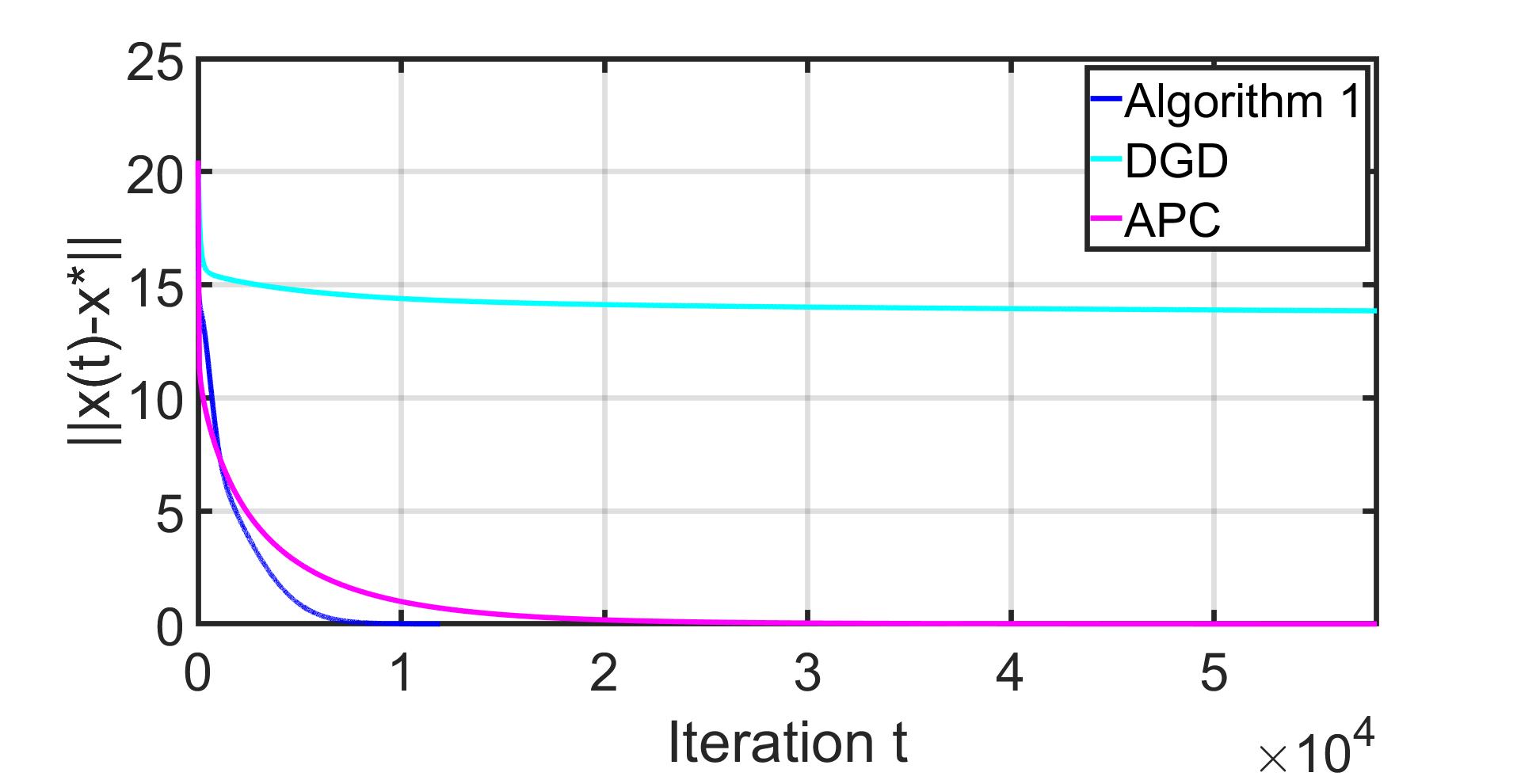}
  \caption{}
  \label{fig:comp_opt}
\end{subfigure}
\caption{\footnotesize{Temporal evolution of error norm for estimate $\norm{x(t)-x^*}$ against number of iterations for Example 1, under Algorithm 1, DGD and APC~\cite{azizan2019distributed}; with (a) arbitrary parameter choices (b) optimal parameter choices. Initialization for (a) and (b) both: (Algorithm 1) $x(0) = [0,\ldots,0]^T$, $K(-1) = O_{n \times n}$; (DGD) $x(0) = [0,\ldots,0]^T$; (APC) according to the algorithm. In (a): (Algorithm 1) $\beta = 100, \alpha = 2*10^{-7}, \delta = 1$; (DGD) $\delta = 10^{-7}$; (APC) $\gamma = \eta = 1$. In (b): (Algorithm 1) $\beta = 5, \alpha = 3.17*10^{-7}, \delta = 1.95$; (DGD) $\delta = 3.17*10^{-7}$; (APC) $\gamma^*=1.08, \eta^*=12.03$.}}
\label{fig:comp}
\end{figure*}

\subsection{Comparison with DGD}

Now we are ready to present our key result, which is a formal comparison between the convergence speed of Algorithm~\ref{algo_1} and the original DGD method in~\eqref{eqn:dgd}. We define a few necessary notations.
\begin{itemize}
\setlength\itemsep{0.5em}
    \item Let, the estimate of the optimum point computed by Algorithm 1 and DGD be denoted by $x_1(t)$ and $x_2(t)$, respectively, after $t$ iterations.
    \item Let, $z_i(t)=x_i(t)-x^*$ for $i=1,2$.
    \item Let, the known upper bound on $\norm{z_i(t)}$ be $E_i(t)$.
\end{itemize}
~

\noindent \fbox{\begin{minipage}{0.47\textwidth}
\begin{theorem} \label{thm:thm1}
Consider Algorithm 1 and the DGD algorithm~\eqref{eqn:dgd} with identical initial estimate $x(0) \in \R^n$. If Assumption 1 holds, then there exists $t_{sw} < \infty$ such that $E_1(t) < E_2(t)$ for all $t> t_{sw}$.
\end{theorem}
\end{minipage}}
~\\


Theorem~\ref{thm:thm1} implies that if the server executes both Algorithm 1 and the original DGD algorithm using the same initial estimate $x(0)$ of $x^*$, then after a certain number of iterations the upper bound on the error-norm generated by Algorithm 1 is smaller than the original DGD. \\

\begin{table*}[htb!]
\caption{}
\begin{center}
\begin{tabular}{|p{3cm}|p{1cm}|p{1.5cm}|p{5cm}|}
\hline
\multicolumn{4}{|c|}{Comparison Between DGD, APC and Algorithm~\ref{algo_1} with Optimal Parameters on \textit{``bcsstm07''};} \\
\multicolumn{4}{|c|}{regarding (a) Number of Iterations to Attain a Relative Estimation Error $10^{-4}$,} \\
\multicolumn{4}{|c|}{ and (b) Decay rate of the Instantaneous Error-norm. Here, $cond(\A) = 5.8 * 10^{7}$.} \\
\hline
\hline
 Algorithm & DGD & APC & Algo.~\ref{algo_1} ($\beta = 5, \, K(-1) = O_{n \times n}$) \\
\hline
(a) Iterations needed & $> 10^5$ & $4.85*10^4$ & $2.11 * 10^4$ \\
\hline
(b) Rate of Decrease & $0.9999$ & $0.9672$ & $0.9583 + 4.3 * 10^8 * (0.9999)^{t+1}$ \\
\hline
\end{tabular}
\end{center}
\label{tab:comp}
\end{table*}

\section{EXPERIMENTS}
\label{sec:exp}

In this section, we present the experimental results to verify the obtained theoretical convergence guarantees of Algorithm~\ref{algo_1}. The matrix $A$ is loaded from one of the real data-sets in SuiteSparse Matrix Collection\footnote{https://sparse.tamu.edu}. Our primary focus here is obtaining faster convergence on a problem~\eqref{eqn:opt_1} where the condition number of $\A$ is large enough, specifically $5.8 * 10^{7}$, while keeping the computational load minimum.\\

\textit{Example $1$.} The matrix $A$ is a real symmetric positive definite matrix from the benchmark dataset \textit{``bcsstm07''} which is part of a structural problem. The dimension of $A$ is $420 \times 420$. We generate the vector $b$ by setting $b=Ax^*$ where $x^*$ is $n=420$ dimensional vector of one's. Since $A$ is positive definite,~\eqref{eqn:opt_1} has a unique solution $x^*$. The data $(A,b)$ is split among $m=10$ machines, so that $n^i=42$, $i=1,\ldots,10$.\\

We apply Algorithm 1 to solve the aforementioned optimization problem. It can be seen that, the speed of the algorithm depends on the parametric choices of $\beta, \alpha, \delta$ (ref. Fig.~\ref{fig:bcsstm_1a}). Also, the algorithm converges to $x^*$ irrespective of the initial choice of the entries in $x(0)$ and $K(-1)$ (ref. Fig.~\ref{fig:bcsstm_1b}).\\

The experimental results have been compared with the conventional distributed gradient descent (DGD) method and the accelerated projection-based consensus (APC)~\cite{azizan2019distributed} method (ref. Fig.~\ref{fig:comp}). The pre-conditioning technique for the distributed heavy-ball method in~\cite{azizan2019distributed} has the same theoretical rate as APC. Besides, it has an additional computational overhead of $\bigO(p^2 n + p^3)$, with $p \leq n_i$ depending on the given matrix $A_i$, which is larger than the cost of Algorithm~\ref{algo_1} in general. \\

We compare the number of iterations by DGD, APC and Algorithm~\ref{algo_1} to reach a relative estimation error (defined as $\dfrac{\norm{x(t)-x^*}}{\norm{x^*}}$) of $10^{-4}$ (ref. Table~\ref{tab:comp}). The tuning parameter $\beta$ has been set at $5$ for Algorithm 1. The rest of the algorithm parameters have been set such that the respective algorithms will have their smallest possible convergence rates. Specifically, $(\alpha^*=3.17*10^{-7}, \delta^*=1.95)$ for Algorithm 1, $(\gamma^*=1.08, \eta^*=12.03)$ for APC and $\delta^*=3.17*10^{-7}$ for DGD. We found that Algorithm~\ref{algo_1} performs the fastest among these algorithms. Note that evaluating the optimal tuning parameters for any of these algorithms requires knowledge about the smallest and largest Eigenvalues of $\A$. \\

We notice that, the error norm $\norm{x(t)-x^*}$ for Algorithm~\ref{algo_1} is less than that of DGD from an approximate iteration index of $300$ onward (ref. Fig.~\ref{fig:comp_opt}). Since $x(0)$ is same for both of the algorithms, this observation is in agreement with the claim in Theorem~\ref{thm:thm1}. Similarly, Algorithm~\ref{algo_1} is faster than APC after $1114$ iterations (ref. Fig.~\ref{fig:comp_opt}).

\section{SUMMARY}

In this paper, we have proposed an algorithm for solving distributed linear least-squares minimization problems over a federated architecture with minimum computational load. However, the algorithm can be emulated in a rooted peer-to-peer network. The vital contribution lies in mitigating the detrimental impact of ill-conditioning on the convergence of the traditional gradient descent method. The computation of the pre-conditioner is done at the server level without requiring any access to the data. In practice, we test the algorithm on a real data-set with a significant condition number ($cond(\A) = 5.8 * 10^{7}$) and get much better performance compared to the classical distributed gradient algorithm as well as the recently proposed accelerated projection-consensus algorithm (APC), regarding the number of iterations needed for convergence to the true solution. We have formally shown the proposed algorithm to converge faster than the distributed gradient method, while the APC method only speculates it to be faster.

\section*{Acknowledgements}
This work is being carried out as a part of the Pipeline System Integrity Management Project, which is supported by the Petroleum Institute, Khalifa University of Science and Technology, Abu Dhabi, UAE. Nirupam Gupta was sponsored by the Army Research Laboratory under Cooperative Agreement W911NF- 17-2-0196. 

\bibliographystyle{unsrt}
\bibliography{refs}

\begin{thebibliography}{10}

\bibitem{lynch1996distributed}
Nancy~A Lynch.
\newblock {\em Distributed algorithms}.
\newblock Elsevier, 1996.

\bibitem{tanenbaum1981network}
Andrew~S Tanenbaum.
\newblock Network protocols.
\newblock {\em ACM Computing Surveys (CSUR)}, 13(4):453--489, 1981.

\bibitem{yang2019federated}
Qiang Yang, Yang Liu, Tianjian Chen, and Yongxin Tong.
\newblock Federated machine learning: Concept and applications.
\newblock {\em ACM Transactions on Intelligent Systems and Technology (TIST)},
  10(2):1--19, 2019.

\bibitem{bertsekas1989parallel}
Dimitri~P Bertsekas and John~N Tsitsiklis.
\newblock {\em Parallel and distributed computation: numerical methods},
  volume~23.
\newblock Prentice hall Englewood Cliffs, NJ, 1989.

\bibitem{azizan2019distributed}
Navid Azizan-Ruhi, Farshad Lahouti, Amir~Salman Avestimehr, and Babak Hassibi.
\newblock Distributed solution of large-scale linear systems via accelerated
  projection-based consensus.
\newblock {\em IEEE Transactions on Signal Processing}, 67(14):3806--3817,
  2019.

\bibitem{nesterov27method}
Y~Nesterov.
\newblock A method of solving a convex programming problem with convergence
  rate {$O$}($1/k^2$).
\newblock {\em Sov. Math. Doklady}, 27(2):372--376, 1983.

\bibitem{fessler2008image}
Jeffrey~A Fessler.
\newblock Image reconstruction: Algorithms and analysis.
\newblock \url{http://web.eecs.umich.edu/~fessler/book/c-opt.pdf}.
\newblock [Online book draft; accessed 17-February-2020].

\bibitem{meijerink1977iterative}
J~Andvandervorst Meijerink and Henk~A Van Der~Vorst.
\newblock An iterative solution method for linear systems of which the
  coefficient matrix is a symmetric {$M$}-matrix.
\newblock {\em Mathematics of computation}, 31(137):148--162, 1977.

\bibitem{axelsson1985survey}
Owe Axelsson.
\newblock A survey of preconditioned iterative methods for linear systems of
  algebraic equations.
\newblock {\em BIT Numerical Mathematics}, 25(1):165--187, 1985.

\bibitem{benzi2002preconditioning}
Michele Benzi.
\newblock Preconditioning techniques for large linear systems: a survey.
\newblock {\em Journal of computational Physics}, 182(2):418--477, 2002.

\end{thebibliography}

\section*{APPENDIX}

\subsection{Proof of Lemma 1}

From~\eqref{eqn:Rij} and noticing that $\sum_{i=1}^m R^i_j(t-1) = \left(\A+\beta I\right)k_j(t-1) - e_j$,
dynamics~\eqref{eqn:kcol_update} can be rewritten as
\begin{align} \label{eqn:kcol_2}
    k_j(t) & = k_j(t-1) - \alpha \left[\left(\A+\beta I\right)k_j(t-1) - e_j\right].
\end{align}
Define $\Tilde{k}_j(t) := k_j(t) - k_j^*$. From the definition of $K^*$ we have,
\begin{align*}
    \left(\A+\beta I\right)K^* = I \implies \left(\A+\beta I\right) k_j^* = e_j, \, j=1,...,n.
\end{align*}
From~\eqref{eqn:kcol_2},
\begin{align}
    \Tilde{k}_j(t) = \left[I- \alpha \left(\A+\beta I\right)\right] \Tilde{k}_j(t-1).
\end{align}
Since $(\A+\beta I)$ is positive definite for $\beta >0$, $\exists \alpha$ for which there is a positive $\rho_K < 1$ such that $\norm{\Tilde{k}_j(t)} \leq \rho_K \norm{\Tilde{k}_j(t-1)}$, $t=0,1,2,...$, where the smallest value of $\rho_K$ is $\dfrac{\kappa(\A+\beta I)-1}{\kappa(\A+\beta I)+1}$ (ref. Corollary 11.3.3 and Chapter 11.3.3 of~\cite{fessler2008image}). As $\kappa(\A+\beta I)=\dfrac{\lambda+\beta}{\gamma+\beta}$, the claim follows.

\subsection{Proof of Lemma 2}

Consider the Eigen decomposition $\A = U \Sigma U^T$, where $U$ is the Eigenvector matrix of $\A$ and $\Sigma$ is a diagonal matrix with the Eigenvalues of $\A$ in the diagonal. Since $K^* := \left(\A+\beta I\right)^{-1}$ and the identity matrix can be written as $I=U U^T$, we have 
\begin{align*}
    & K^* = U diag\{\dfrac{1}{\lambda_{k}+\beta}\}_{k=1}^n U^T \\
    \implies & K^*\A = U diag\{\dfrac{1}{\lambda_{k}+\beta}\}_{k=1}^n U^T U diag\{\lambda_{k}\}_{k=1}^n U^T \\
    \implies & K^*\A = U diag\{\dfrac{\lambda_{k}}{\lambda_{k}+\beta}\}_{k=1}^n U^T,
\end{align*}
where $\lambda_{k} > 0$ are the Eigenvalues of the positive definite matrix $\A$. Since $\dfrac{d}{d\lambda_k} \left(\dfrac{\lambda_{k}}{\lambda_{k}+\beta}\right) > 0$, the claim follows.

\subsection{Proof of Proposition 1}

Using~\eqref{eqn:g_i} and observing that $A^T A = \sum_{i=1}^m (A^i)^T A^i$, dynamics~\eqref{eqn:x_update} can be rewritten as
\begin{align}
    x(t+1) & = x(t) - \delta K(t) A^T\left(Ax(t)-b\right). \label{eqn:x_central}
\end{align}
Define $z(t):= x(t) - x^*$. The objective cost in~\eqref{eqn:opt_1} can be rewritten as $\dfrac{1}{2}\norm{A \, x -b}^2$, gradient of which is given by $\A x - A^T b$. Thus, $x^*$ satisfies $\A \, x^* = A^T b$.
From~\eqref{eqn:x_central} and $\A \, x^* = A^T b$, we get
\begin{align}
    z(t+1) & = z(t) - \delta K(t) \A (x(t) - x^*) \nonumber \\
    & = \left(I- \delta K(t)\A\right) z(t) \nonumber \\
    & = \left(I- \delta K^*\A\right)z(t) - \delta \Tilde{K}(t)\A z(t) \nonumber \\
    & = \left(I- \delta K^*\A\right)z(t) - \delta u(t),  \label{eqn:z}
\end{align}
where $\Tilde{K}(t)$ comprises of the columns $\Tilde{k}_j(t) := k_j(t) - k_j^*$, $j=1,...,n$ and $u(t) := \Tilde{K}(t)\A z(t)$.
From Lemma~\ref{thm:lem1}, it follows that
\begin{align}
    & \norm{\Tilde{k}^j(t)}^2 \leq (\rho^*_K)^{2t+2} \norm{\Tilde{k}^j(-1)}^2 \nonumber \\
    \implies & \norm{\Tilde{K}(t)}_F^2 \leq (\rho^*_K)^{2t+2} \norm{\Tilde{K}(-1)}_F^2. \label{eqn:K_frob}
\end{align}
Now, 
\begin{align}
    \norm{u(t)} & \leq \norm{\Tilde{K}(t)} \norm{\A} \norm{z(t)} \nonumber \\
    & \leq \norm{\Tilde{K}(t)}_F \norm{\A} \norm{z(t)} \nonumber \\
    & \leq \lambda \norm{\Tilde{K}(-1)}_F (\rho^*_K)^{t+1} \norm{z(t)}, \label{eqn:u_norm}
\end{align}
where the last inequality follows from~\eqref{eqn:K_frob}. From~\eqref{eqn:z} and~\eqref{eqn:u_norm},
\begin{align}
    \norm{z(t+1)} & \leq \left(\norm{I- \delta K^*\A} + \sigma_0 (\rho^*_K)^{t+1}\right) \norm{z(t)}, \, \forall t. \label{eqn:z_norm}
\end{align}
Since $K^*\A$ is positive definite for $\beta >0$ (from Lemma~\ref{thm:lem2}), $\exists \delta$ for which $\rho_{\beta} := \norm{I- \delta K^*\A} < 1$ and~\eqref{eqn:z_norm} holds. The smallest value of $\rho_{\beta}$ is $\dfrac{\kappa(K^*\A)-1}{\kappa(K^*\A)+1}$. Then with the largest and the smallest Eigenvalues of $K^*\A$ obtained in Lemma~\ref{thm:lem2}, the claim follows.

\subsection{Proof of Theorem 1}

Define, $s_t := \left(\ro+\sigma_0 (\rho^*_K)^{t+1}\right)$. It can be easily checked that $\ro < \rho_{GD}$ for $\beta > 0$. Since $\rho_K^* < 1$ and $\ro < \rho_{GD}$, there exists $\tau < \infty$ such that $s_t < \rho_{GD}, \, \forall t > \tau$.
Define, $r_t := s_t/\rho_{GD}$. Then $r_t < 1 \, \forall t > \tau$.
From the recursion in Proposition~\ref{thm:pro1}, we have
\begin{align}
   \norm{z_1(t+1)} & \leq \left(\Pi_{k=\tau+1}^t \, s_k\right) \norm{z_1(\tau + 1)}, \, \forall t > \tau. \label{eqn:z1}
\end{align}
Since $\{s_t > 0\}_{t \geq 0}$ is a strictly decreasing sequence, from the same recursion we also have
\begin{align}
    & \norm{z_1(\tau+1)} \leq s_{\tau} \norm{z_1(\tau)} \leq s_{0} \norm{z_1(\tau)} \nonumber \\
    \implies & \norm{z_1(\tau+1)} \leq s_{0}^{\tau+1} \norm{z(0)}. \label{eqn:z2}
\end{align}
Combining~\eqref{eqn:z1} and~\eqref{eqn:z2},
\begin{align*}
    \norm{z_1(\tau+1)} & \leq \left(\Pi_{k=\tau+1}^t \, s_k\right)s_{0}^{\tau+1} \norm{z(0)} \leq (s_{\tau +1})^{t-\tau} s_{0}^{\tau+1} \norm{z(0)}\\
    & = (r_{\tau +1} \rho_{GD})^{t-\tau} s_{0}^{\tau+1} \norm{z(0)} \\
    & = c_1(r_{\tau+1} \rho_{GD})^{t+1}\norm{z(0)}, \, \forall t > \tau
\end{align*}
where $c_1 = \left(\dfrac{s_{0}}{r_{\tau +1} \rho_{GD}} \right)^{\tau+1}$ is a constant.
Thus, we have
\begin{align*}
    E_1(t+1) & = c_1(r_{\tau+1} \rho_{GD})^{t+1}\norm{z(0)} \\
    & = c_1(r_{\tau+1})^{t+1} E_2(t+1), \, \forall t > \tau.
\end{align*}
Since $r_{\tau+1} < 1$ and $c_1$ is a constant, $\exists \, t_{sw} < \infty$ such that $c_1(r_{\tau+1})^{t+1} < 1 \, \forall t>t_{sw}$.
This completes the proof.


\addtolength{\textheight}{-12cm}

\end{document}